\documentclass[opre,nonblindrev]{informs3} 

\DoubleSpacedXI 

\usepackage{endnotes}
\let\footnote=\endnote

\usepackage[table,xcdraw]{xcolor}
\usepackage{longtable}

\usepackage{natbib}
 \bibpunct[, ]{(}{)}{,}{a}{}{,}%
 \def\bibfont{\small}%
\TheoremsNumberedThrough     
\ECRepeatTheorems

\EquationsNumberedThrough    

\MANUSCRIPTNO{2} 

\begin{document}


\RUNAUTHOR{Lazzaro, Foraker, and Curry}

\RUNTITLE{Brigade Reassignment Problem}

\TITLE{The Service Academy Brigade Reassignment Problem}

\ARTICLEAUTHORS{%
\AUTHOR{Gary Lazzaro, Jay Foraker, Rob Curry}
\AFF{Department of Mathematics, United States Naval Academy, Annapolis, MD, 21402, \\ \EMAIL{lazzaro@usna.edu,}\EMAIL{foraker@usna.edu,}\EMAIL{rcurry@usna.edu}} 

} 

\ABSTRACT{%
After the COVID pandemic, student leadership development within the brigade of midshipmen at the United States Naval Academy (USNA) was severely degraded due to reduced student interaction from isolation and remote learning.  Brigade leadership decided that reassignment of some students from their previous companies to new companies could facilitate new leadership bonds for the next academic year.  All students in the classes of 2023 and 2024 are reassigned from their previous companies to new companies with the Brigade Reassignment Problem (BRP). The first goal of BRP is to minimize the number of students reassigned to their current company.  The second goal seeks to improve the homogeneity of each company in terms of average metrics defined by brigade leadership. We create three similar mathematical programming models to achieve different goals of the BRP. The first model focuses on the first goal and quickly obtains an optimal solution. A second model focuses on the second goal and obtains high-quality feasible solutions. A third model combines the two goals and obtains an optimal solution within seven hours.  USNA leadership implemented results from the first model to reassign the classes of 2023 and 2024 to new companies for the next academic year.\\
\textit{Composite Group:} Resource, Readiness and Training\\
\textit{Military OR Application Area:} Analytic Support to Training and Education\\
\textit{OR Methodology:} Linear/Integer Programming
}%

\HISTORY{Accepted 11 Oct 2021, Revision 1 submission, 03 October 2021, First Submission, 29 June 2021.}

\maketitle
\section{INTRODUCTION}
\label{S:Intro}

\begin{figure}[ht]
\centering
\includegraphics[scale=0.35]{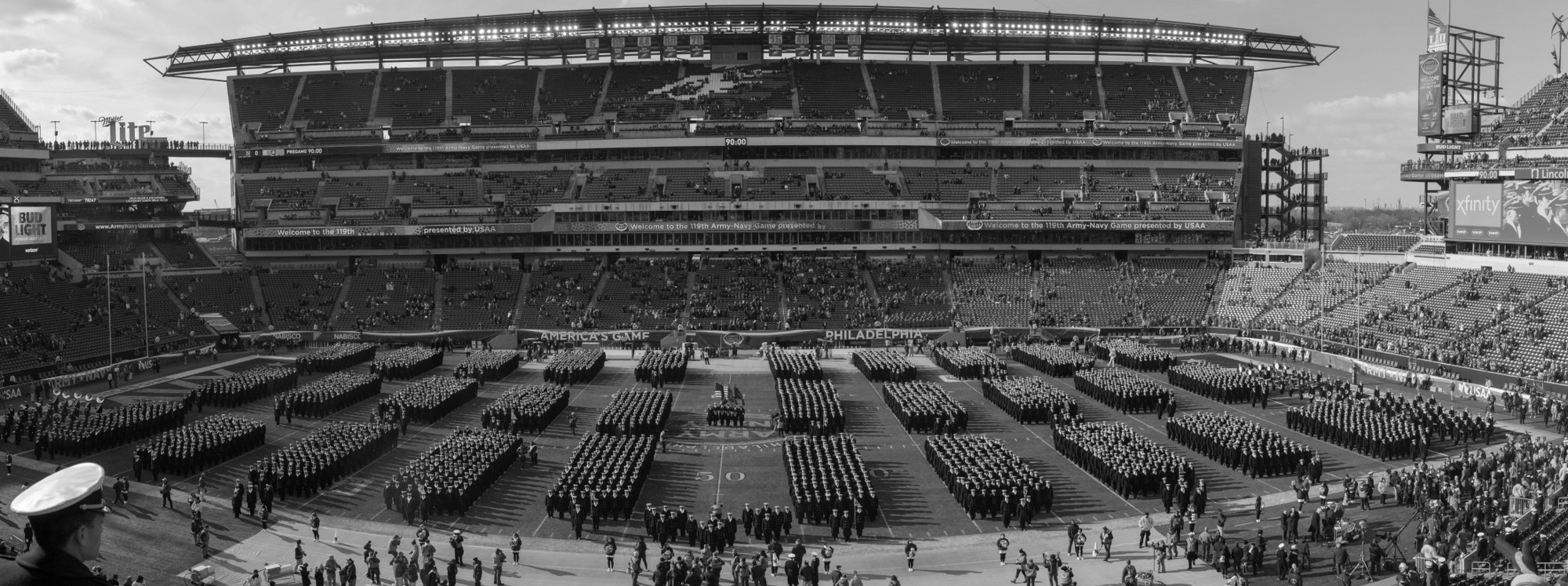}
\caption{Brigade of Midshipmen in Formation \citep{DVIDS:2018}}
\label{fig:brigadeformation}
\end{figure}

The United States Naval Academy (USNA) has approximately 4,400 students referred to as midshipmen.  As of spring 2021, there are 1,165 freshmen, 1,097 sophomores, 1,130 juniors, and 1,050 seniors. The students are arranged into a midshipmen brigade of 6 battalions which are further subdivided into 5 companies each for a total of 30 companies. Each company has between 33-42 students from each class.  Figure \ref{fig:brigadeformation} shows the Brigade of Midshipmen in formation of companies on a field. The brigade organization exists to develop leadership traits with different levels of responsibility throughout the four years a student attends USNA \citep{atwaterAndYummarino1989}.  


During the COVID-19 pandemic, midshipmen were isolated from their leadership experiences within the brigade.  Restrictions due to COVID forced all midshipmen home for remote learning during Spring 2020 for full virtual learning.  Summer training at sea was cancelled in 2020, and all students completed virtual summer school. The brigade returned in Fall 2020, but various isolation protocols and hybrid classroom learning continued through the end of spring semester 2021.  Furthermore, USNA leadership felt the pressures from the COVID environment exhibited strains to individual and team morals analogous to mechanical brittle fractures in material science \citep{Sears:2020}. Specifically, COVID protocols may have inadvertently created undetected stress or tension between individual students or groups that could result in cracks or failures in leadership relationships. USNA leaders believe company-wide reassignment of students could restore a healthy environment for leadership development.

Senior leadership decided to reassign those classes of midshipmen most adversely affected by the COVID protocols to new companies in order to reset the brigade to normal operations for the next academic year. We define the problem of determining to which new company each midshipman should be reassigned as the Brigade Reassignment Problem (BRP). The class of 2021 students (seniors) are not included in the BRP since they will not participate in the next academic year because of graduation and commissioning. The class of 2022 students (juniors) are also not included in the BRP since they will remain in their current companies in order to form the student leadership cadre for the following academic year. Senior leaders believe these students did not have an opportunity to develop any significant bonds with the underclassmen they will lead in the next academic year. Students in the classes of 2023 (sophomores) and 2024 (freshmen) lost many leadership development opportunities due to virtual learning and COVID isolation. Thus, the BRP focuses solely on reassigning the classes of 2023 and 2024. The reassignment problem is colloquially referred to as the brigade `shotgun' as 2,262 midshipmen are scattered to new companies throughout the brigade.  

\section{PROBLEM DESCRIPTION} \label{prob-description}

\begin{figure}[ht]
\centering
\includegraphics[scale=0.95]{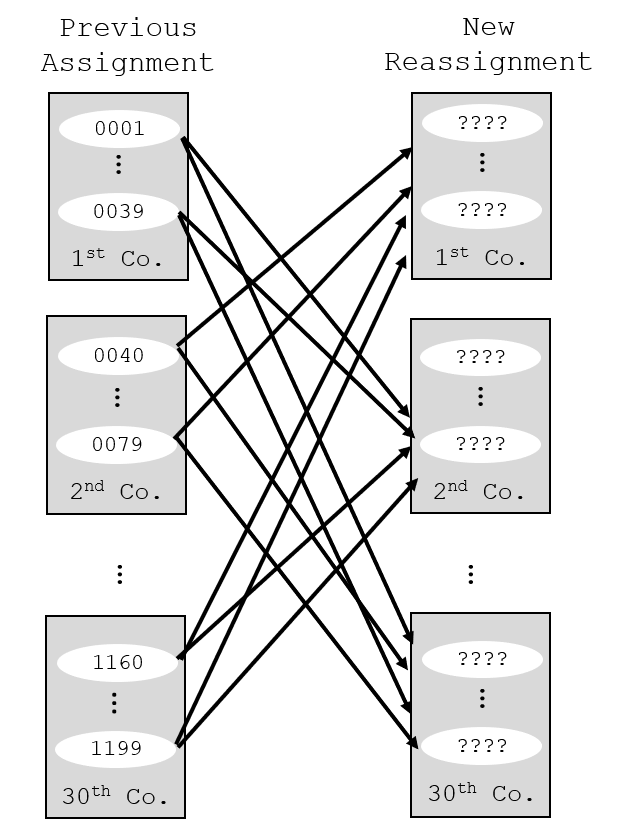}
\caption{Visualization of the Brigade Reassignment Problem}
\label{fig:reassignmentpic}
\end{figure}

The BRP is an enhanced modification of the classic assignment problem \citep{kuhn1955}. The BRP focuses on reassigning one year of students at a time, which is 25\% of the brigade.  One class year of students begins with assignments of each student to the 30 different companies within the brigade. Each company consists of 33 to 42 students from each class.  The existing variability in company size is due to attrition of students over time. The BRP requires all students in a designated class to be reassigned to new companies.  Figure \ref{fig:reassignmentpic} shows each student as a 4-digit number assigned to a company on the left side.  Each student will need to be reassigned to a different company on the right side, preferably a different numbered company than the previous assignment. A simple but naive way to perform such a reassignment would be to randomly pick new companies for each student, but the BRP desires to obtain a better result than random reassignment.

Solving the BRP via mathematical programming seeks to achieve two goals. The primary goal of BRP is to minimize the number of students that remain in the same numbered company after reassignment. Ideally, no student will be reassigned to the same company as before.  
The secondary goal of BRP is to improve the homogeneity of each company in the brigade defined as the differences in average Academic Order of Merit (AOM) and Military Order of Merit (MOM) scores between companies.  USNA leadership desires balanced companies in order to avoid the the perception that any particular company is known to consist of only `smart', `squared-away', or `struggling' students. Each company has upper and lower limits on the average scores of students for academics and military bearing. 

In addition to these two goals, the BRP must observe a variety of other requirements on company composition. Each company possesses requirements on the upper and lower limits on the total number of midshipmen, average physical fitness scores, the number of students assigned to a specific task force, and the number of students with prior military service. Task force students are identified to indoctrinate the new class reporting in the summer. Students having prior military service are those who were enlisted in the military before enrolling at USNA.  Leadership also requires some improvement to the average diversity metrics of companies, so we include restrictions that seek to balance companies with respect to gender and race. Leadership provided a maximum limit on the number of student athletes from any of the large varsity or club sports in each company.  In response to the recent focus on sexual assault prevention and response (SAPR), the BRP ensures that each new company is assigned a minimum number of students trained as a SAPR guide.  Next, there are exactly 30 international students, and the BRP ensures exactly one international student is in each company. To avoid conflict of interest in military activities, the BRP ensures that students involved in a romantic dating relationship are not assigned to the same company. Finally, a small group of students are under specific supervision from a battalion officer, and these students must be assigned to a different company within the same battalion.

\section{PREVIOUS WORK}

The BRP is a unique variation and application of the Generalized Assignment Problem (GAP) which involves minimizing the total cost of assigning a set of agents to a set of tasks where each agent/task combination is assigned a specific cost. Each agent could be assigned to multiple tasks. A common GAP objective function maximizes the total profit of assigning agents to tasks subject to budgetary constraints on each agent. The GAP is known to be NP-hard, and is often modeled and solved as an integer program (IP) \citep{nauss2003}. A common variation of the GAP gives each agent a budget and assignment cost of one unit, and efficient algorithms exist for optimal solutions \citep{kuhn1955}. Other variations and extensions to the GAP include the quadratic assignment problem \citep{Lawler1963}, the transportation problem \citep{FordAndFulkerson1956}, and the dynamic assignment problem \citep{SpiveyAndPowell2004}.


The BRP can be viewed as a variation of the GAP in which multiple midshipmen are assigned to each company. However, the BRP also contains a non-trivial set of side constraints, and the objective is not as straightforward as minimizing the total cost of assignments. Instead, our models of BRP focus on different ways to approach an equitable assignment plan subject to a unique set of side constraints. The notion of equity varies for each model.

We briefly discuss defense-focused applications of the GAP related to the BRP. The most common military application of the GAP is the weapon-target assignment problem, in which weapons are assigned to a set of adversarial targets subject to capability constraints. Some works solve this problem heuristically using genetic algorithms \citep{LeeEtAl2003}, while others solve this problem exactly \citep{AhujaEtAl2007}. More recent works linearize this problem to present an alternative exact approach \citep{LuAndChen2021}. Similar to the BRP, the weapon-target assignment problem also contains a set of complicating side constraints that make it unique when compared to the GAP. The BRP objective differs from the weapon-target assignment problem in that no straightforward cost exists for each midshipman/company combination.  

\cite{LiangAndBuclatin1988} apply similar techniques to the BRP for solving Navy training resource problems. The authors develop a network-base model similar to Figure \ref{fig:reassignmentpic} that assigns personnel to their assignments with en route training. Their network model also includes a set of complicating constraints that model training school class capacities. This model focuses on overall staff and training utilization while also considering assignment costs. Unlike the BRP, their work does not seek to balance the capabilities and/or competencies of the staff assigned to each training resource. 

A similar line of research to the BRP studies personnel assignment problems in military settings. \cite{SonmezAndSwitzer2013} study matching branch-of-choice contracts with student cadets at the United States Military Academy via a market design approach. More generally, \cite{gass1991} applies personnel decision-making to military planning with a high-level description of Markov models, network models, and goal-programming models. Furthermore, the author shows how these models can be used in military personnel planning applications. \cite{rigopoulosEtAl2011} solve a military staff assignment problem, in which the authors propose an approach that utilizes multi-criteria analysis to assign staff according to various predefined criteria instead of a single objective function. 

\cite{TorosluAndArslanoglu2007} study multi-objective personnel assignment outside of a military context. They propose a variation of the GAP having side constraints that require some agents to be assigned to only some strict subset of possible tasks. Rather than considering each objective independently (as in our work), these authors take a multi-objective approach to both minimize the number of violated side constraints and maximize the summation of the weighted assignment. Their solution procedure involves a genetic algorithm heuristic as opposed to our approach of solving each different IP to optimality or near-optimality. 

Other non-defense application areas of the GAP include manufacturing \citep{LeeAndKim1998}, scheduling \citep{Carreno1990}, education \citep{FaudziEtAl2018}, and the peer-review process \citep{HartvigsenEtAl1999}. For a full analysis and survey of the GAP literature, see \cite{oncan2007}.  


\section{MODEL FORMULATION}

We describe our mathematical programming formulation for creating different versions of the BRP. Without the presence of additional side restrictions, a class of 1,100 students can be distributed amongst the 30 companies of the brigade of midshipmen in $30^{1100}$ different ways. Due to this complexity and the combinatorial number of reassignment plans, USNA leadership  previously created a company reassignment plan by hand, without the tools of mathematical programming. The reassignment process took days to complete. As long as the tolerances for each of the requirements described in Section \ref{prob-description} are set to reasonable values, many feasible reassignments exist. A naive model without an objective function could produce a sub-optimal solution to reassign all students to new companies subject to these requirements.

We create three unique optimization models for the BRP to achieve different goals in order to allow brigade leadership to achieve broader goals for the composition of the companies within the brigade.  All three models have many common sets, parameters, decision variables, and constraints.  Each model has a unique objective function designed to achieve a specific goal.  Our first model, referred to as \textit{BRP-MIN} is designed to achieve the first goal of the reassignment, which is to minimize the number of midshipmen reassigned to their previous company.  \textit{BRP-MIN} is simple and quickly produces an optimal solution.  The second and third models are designed to meet the secondary goals of reassignment, which improves the homogeneity between companies.  Our second model, referred to as \textit{BRP-DEV}, seeks to minimize the absolute deviation of company averages for academic and military order of merit scores.  The third model, referred to as \textit{BRP-PAIRS} seeks to minimize the number of student pairs reassigned to a new company that were also paired in the same previous company.  

\subsection{Sets and Indices for all Models}

\begin{tabbing}
\hspace*{.3cm} \= $ a \in A $ \hspace{2.0cm} \= set of student alpha codes for a single class of students.\\
\> $ g \in G $ \> set of genders, $G=\{male, female\}$.\\
\> $ e \in E $ \> set of Navy defined race/ethnicity codes \citep{NPC:2021}. \\
\> $ c \in C $ \> set of companies, $C=\{C1, C2, ...,C30\}$.\\
\> $ b \in B $ \> set of battalions, $B=\{B1, B2,...,B6\}$.\\
\> $ w \in W_b $ \> subset of five companies within battalion $b \in B$, $W_b \subset C$.\\
\>\> $W_1 = \{C1, C2, C3, C4, C5\}$ , $W_2 = \{C6, C7, C8, C9, C10\}$, etc. \\
\> $ q \in Q $ \> set of student attributes to be counted in each company, \\
\>\> $Q=$\{\textit{all students, task force students, prior service students}\}.\\
\> $ m \in M $ \> set of merit scores for a student, \\
\>\> $M=$\{\textit{AOM score, MOM score, PRT score}\}.\\
\> $ v \in V $ \> set of 13 varsity or club sports with a large number of students,\\
\>\>$V=$\{\textit{heavyweight crew, lightweight crew, women's crew, football, sprint football,} \\
\>\> \textit{men's lacrosse, women's lacrosse, offshore sailing, men's rugby, women's rugby},\\
\>\> \textit{swimming \& diving, track, wrestling}\}.\\
\> $s \in S_v$ \> subset of student alpha codes that participate in sport $v$, $S_v \subset A \quad \forall v \in V$. \\
\>$d \in D_q$ \> subset of student alpha codes that are designated members of group $q \in Q$, $D_q \subseteq A$.\\
\> $p \in P$ \> subset of student SAPR guide alpha codes, $P \subset A$.\\
\> $ pn \in PN$ \> subset of companies requiring a SAPR guide, $PN \subset C$.\\
\> $ is \in IS $ \> subset of international student alpha codes, $IS \subset A$.\\
\> $ in \in IN$ \> subset of companies requiring an international student, $IN \subset C$.\\
\> $ k \in K $ \> subset of student alpha codes that must remain within the same battalion, $K \subset A$.\\
\> $ (r,\rho) \in R$ \> subset of student alpha code pairs that must be assigned to different companies.
\end{tabbing}

\subsection{Parameters for all Models}

Each student has a broad range of qualities that should be found in all of the companies of the brigade.  Brigade leadership desires to have companies with a quality spread of students having these characteristics.  USNA leadership defines minimum and maximum tolerances for many student attributes within any particular company as well.     

\begin{tabbing}
\hspace*{.3cm} \= $score_{a,m}$ \hspace{2.0cm} \=  merit score for student $a \in A$ in category $m \in M$. \\
\> $gender\_demo_{a,g}$ \>1, if student $a \in A$ with gender $g \in G$ is male, and 0, if female.\\
\> $race\_demo_{a,e}$ \>1, if student $a \in A$ with ethnicity $e \in E$ is white, and 0, for all other groups.\\
\> $old\_comp_{a,c}$ \> 1, if student $a \in A$ was previously a member of company $c \in C$, and 0, otherwise.\\
\> $min\_sapr$ \> minimum number of SAPR guides required in a company.\\
\> $min\_number_{q}$ \> minimum number of students allowed in a company with quality $q \in Q$.\\
\> $max\_number_{q}$ \> maximum number of students allowed in a company with quality $q \in Q$.\\
\> $min\_avg\_score_{m}$ \> minimum allowed average merit score $m \in M$ for a company.\\
\> $max\_avg\_score_{m}$ \> maximum allowed average merit score $m \in M$ for a company.\\
\> $min\_gender_{g}$ \> minimum percentage of students in a company with gender $g \in G$. \\
\> $max\_gender_{g}$ \> maximum percentage of students in a company with gender $g \in G$. \\
\> $min\_race_{e}$ \> minimum percentage of students in a company with race $e \in E$. \\
\> $max\_race_{e}$ \> maximum percentage of students in a company with race $e \in E$. \\
\> $max\_athlete_{v}$ \> maximum number of athletes in a company from sport $v \in V$.
\end{tabbing}

\subsection{Decision Variables for all Models}
We define binary decision variable $x$ to signify to which new company each student is assigned.
\begin{tabbing}
\hspace*{.3cm} \= $x_{a,c} \in \{0,1\}$ \hspace{0.8cm} \= 1, if student with alpha $a$ is assigned to company $c$, and 0, otherwise $\forall \ a \in A, c \in C$.
\end{tabbing}

\subsection{Constraints for all Models}
The following constraints enforce variable relationships, company demographics, and brigade policies. We describe each block of constraints separately to explain the model.
\begin{subequations}
\begin{align}
     &  
\label{onecompany}
\sum_{c \in C} x_{a, c} = 1 \qquad \forall a \in A \\
& \label{numericalmax}
\sum_{d \in D_q} x_{d, c} \le max\_number_q \qquad \forall q \in Q, c \in C\\
&\label{numericalmin}
\sum_{d \in D_q} x_{d, c} \ge min\_number_q \qquad \forall q \in Q, c \in C \\
&\label{meritmax}
\sum_{a \in A} score_{a,m} * x_{a, c} \le max\_avg\_score_m * \sum_{a \in A}x_{a, c} \qquad \forall m \in M, c \in C \\
& \label{meritmin}
\sum_{a \in A} score_{a,m} * x_{a, c} \ge min\_avg\_score_m * \sum_{a \in A}x_{a, c} \qquad \forall m \in M, c \in C \\
& \label{gendermax}
\sum_{a \in A}gender\_demo_{a,g}* x_{a, c} \le max\_gender_g * \sum_{a \in A}x_{a, c}\qquad \forall g \in G, c \in C \\
& \label{gendermin}
\sum_{a \in A}gender\_demo_{a,g}* x_{a, c} \ge min\_gender_g * \sum_{a \in A}x_{a, c} \qquad \forall g \in G, c \in C \\
& \label{racemax}
\sum_{a \in A}race\_demo_{a,e}* x_{a, c} \le max\_race_e * \sum_{a \in A}x_{a, c} \qquad \forall e \in E, c \in C \\
& \label{racemin}
\sum_{a \in A}race\_demo_{a,e}* x_{a, c} \ge min\_race_e * \sum_{a \in A}x_{a, c} \qquad \forall e \in E, c \in C \\
& \label{maxsportconst}
\sum_{s \in S_v} x_{s, c} \le max\_athlete_{v} \qquad \forall v \in V, c \in C \\
& \label{datingconst}
x_{r, c} + x_{\rho,c}\le 1 \qquad \forall (r,\rho) \in R, c \in C \\
& \label{saprconst}
\sum_{p \in P} x_{p, pn} \ge min\_sapr \qquad \forall pn \in PN \\
& \label{intlconst}
\sum_{is \in IS} x_{is, in} = num\_intl \qquad \forall in \in IN \\
& \label{samebattconst}
\sum_{\substack{w \in W_b : \\ arg(old\_comp_{k,w}) \in b}} x_{k, w} \ge 1 \qquad \forall k \in K \\
& \label{binaryconst}
x_{a, c} \in \{0,1\} \qquad \forall a \in A, \ c \in C 
\end{align}
\end{subequations}

Constraint \eqref{onecompany} ensures each student is assigned to exactly one company.  Constraints \eqref{numericalmax} and \eqref{numericalmin} enforce  maximum and minimum  limits for the number of total students, task force students, and prior service students in each company. Constraints \eqref{meritmax} and \eqref{meritmin} create upper and lower bounds for the average academic, military, and physical fitness merit scores of all students assigned to each company.  Constraints \eqref{gendermax}--\eqref{racemin} enforce maximum and minimum limits for the number of white students and male students for each company. Constraints \eqref{maxsportconst} limit the maximum number of student athletes in a particular sport assigned to each company. Constraints \eqref{datingconst} prevent a pair of dating students from being assigned to the same new company. Constraints \eqref{saprconst} ensure that each company has at least the minimum number of SAPR representatives. Constraints \eqref{intlconst} ensures that each company has exactly one international student. Constraints \eqref{samebattconst} ensure that the appropriate group of students remain in their current battalion. Finally, constraints \eqref{binaryconst} enforce binary restrictions on $x$. 

\subsection{Model \textit{BRP-MIN}: Minimize Students Assigned to the Same Company}
\textit{BRP-MIN}'s objective function \eqref{objective1} reflects the goal of minimizing the number of students reassigned to their previous company. The preferred outcome from leadership is to have zero midshipmen reassigned to their previous company.  \textit{BRP-MIN} focuses on the first goal of BRP, with limited focus on the secondary goal through enforcement of constraints. \textit{BRP-MIN} is purposely designed to be simpler than the other models in order to obtain an optimal solution quickly.  
\begin{align}
& \mathrm{Minimize }\sum_{a \in A, c \in C} old\_comp_{a, c} * x_{a,c} \label{objective1}
\end{align}

\subsection{Model \textit{BRP-DEV}: Minimize Absolute Deviation of Average Company Metrics}
\textit{BRP-DEV} attempts to achieve the secondary goal of improving company homogeneity by producing a weighted average of the smallest deviation from the average academic and military order of merit scores between any two pairs of newly assigned companies.  Rather than minimizing the total number of students reassigned to the previous company as in \textit{BRP-MIN}, we include additional constraints that prevent students from being reassigned to their previous company.  \textit{BRP-DEV} seeks to balance the composition of companies, in terms of a weighted average of AOM and MOM scores for each company. Typically, standard deviation is used to measure deviation from the mean, but attempting to minimize standard deviation would introduce non-linearity to the model. Instead, we minimize absolute deviation between weighted average of AOM/MOM scores between each combination of company pairs in order to preserve the linearity of the model.  Additional parameters, decision variables, and constraints must be added to \textit{BRP-DEV} to linearize the absolute value of deviations. For more details, see \citep{Granger:2014}.
\subsubsection{\textit{BRP-DEV} Parameters}
\begin{tabbing}
\hspace*{.3cm} \= $AOM\_weight \in [0,1]$ \hspace{2.0cm} \= the weighted importance of AOM. \\
\> $MOM\_weight \in [0,1]$ \> the weighted importance of MOM. \\
\> \> $AOM\_weight + MOM\_weight = 1$.
\end{tabbing}
\subsubsection{\textit{BRP-DEV} Decision Variables}
\begin{tabbing}
\hspace*{.3cm} \= $y_{c,c'} \ge 0$ \hspace{2.0cm} \= difference in AOM values between companies $c$ and $c'$, $\forall c,   c' \in C: c \neq c'$. \\
\> $z_{c,c'} \ge 0$ \> difference in MOM values between companies $c$ and $c'$, $\forall c,  c' \in C : c \neq c'$.
\end{tabbing}

\subsubsection{\textit{BRP-DEV} Objective Function:}
\begin{align}
\textrm{Minimize   } \  & AOM\_weight * \sum_{\substack{c, c' \in C: \\ c \neq c'}} y_{c,c'} + MOM\_weight *\sum_{\substack{c, c' \in C: \\ c \neq c'}} z_{c,c'} 
\label{objective2}
\end{align}
\subsubsection{\textit{BRP-DEV} Constraints}
\begin{subequations}
\begin{align}
    & \label{nocurr}
\sum_{c \in C} old\_comp_{a,c} * x_{a, c} = 0 \qquad \forall a \in A \\
& \label{positiveAOM}
\sum_{a \in A} score_{a,AOM} * x_{a, c} - \sum_{a \in A} score_{a,AOM} * x_{a, c'}\leq y_{c,c'} \qquad \forall  c, c' \in C: c \neq c' \\
  \label{negativeAOM}
- & \left(\sum_{a \in A} score_{a,AOM} * x_{a, c} - \sum_{a \in A} score_{a,AOM} * x_{a, c'}\right)\leq y_{c,c'} \qquad \forall c, c' \in C: c \neq c' \\
 & \label{positiveMOM}
\sum_{a \in A} score_{a,MOM} * x_{a, c} - \sum_{a \in A} score_{a,MOM} * x_{a, c'}\leq z_{c,c'} \qquad \forall c, c' \in C: c \neq c'\\
 \label{negativeMOM}
- & \left(\sum_{a \in A} score_{a,MOM} * x_{a, c} - \sum_{a \in A} score_{a,MOM} * x_{a, c'}\right)\leq z_{c,c'} \qquad \forall c, c' \in C: c \neq c'\\
& \label{nonnegconst}
y_{c, c'}, \ z_{c, c'} \geq 0 \qquad \forall c, \ c' \in C
\end{align}
\end{subequations}

\textit{BRP-DEV}'s objective function \eqref{objective2} reflects the goal of minimizing the maximum absolute deviation of a weighted combination of average AOM and MOM score differences amongst pairs of companies.  Constraints \eqref{nocurr} are added to the model to ensure that no student is allowed to remain in their previously assigned company.  Constraints \eqref{positiveAOM} and \eqref{negativeAOM} are added in order to linearize the modeling of the average difference of AOM scores between pairs of companies.  Constraints \eqref{positiveMOM} and \eqref{negativeMOM} perform a similar function for the average difference of MOM scores between pairs of companies.  Constraint \eqref{nonnegconst} enforces non-negative decision variables.

\subsection{Model BRP-PAIRS: Minimize Pairwise Reassignment From Old to New Company}

USNA leadership decided to reassign students to new companies in hopes that new leadership bonds will be formed as students interact with new students not in their previous company.  Leadership assumes that new relationships will form because only a small group of students in the same previous company will be reassigned to the same new company. An ideal goal would be for each student to be assigned to a new company having no students from their previous company. However, such a solution is not feasible. Each of the 30 companies contains roughly 39 students from each class. It is impossible to distribute the 39 students from a company by placing a maximum of one student in each of the 29 other companies. The brigade would need to have 40 companies for such a reassignment to be feasible. Therefore, some companies must get assigned at least two students from the same previous company. \textit{BRP-PAIRS} seeks to keep the number of pairs of students with the same reassignment to a minimum. See the Appendix for further discussion of this topic. 

\subsubsection{\textit{BRP-PAIRS} Sets}
\begin{tabbing}
\hspace*{.3cm} \= $ (t_i,t_j) \in T $ \hspace{2.0cm} \= Set of student alpha code pairs assigned to the same previous company.\\
\> \> $T = A \times A :( old\_comp_{a, c} = old\_comp_{a',c} = 1 ):a \neq a'\ \forall c \in C$.
\end{tabbing}
\subsubsection{\textit{BRP-PAIRS} Decision Variables}
\begin{tabbing}
\hspace*{.3cm} \= $ u_{t_i,t_j} \in \{0,1\} $ \hspace{2.0cm} \= 1, if student alpha $t_i$ is assigned to same new company as student alpha $t_j$\\
\> \> 0, otherwise.
\end{tabbing}
\subsubsection{\textit{BRP-PAIRS} Objective Function}
\begin{align}
\textrm{Minimize} \sum_{(t_i, t_j) \in T: i \neq j} u_{t_i,t_j}  
\label{objective3}
\end{align}

\subsubsection{\textit{BRP-PAIRS} Constraints}
\begin{subequations}
\begin{align}
    & \label{midcombo}
 u_{t_i,t_j} \geq x_{t_i,c} + x_{t_j,c} -1 \qquad \forall (t_i,t_j) \in T: i \neq j, c \in C. \\
     & \label{combobinary}
 u_{t_i,t_j} \in \{0 ,1\} \qquad \qquad \qquad \forall (t_i,t_j) \in T: i \neq j. 
\end{align}
\end{subequations}

We present an example to explain set $T$ and decision variable $u$. Suppose students $0001$ and $0002$ were both previously assigned to company $11$, but student 0003 was previously assigned to company 12.  The set $T$ includes the pair of students (0001, 0002), but $T$ does not include the pairs (0001, 0003) and (0002, 0003). If students 0001 and 0002 are both reassigned to the same new company (e.g. company 20), then variable $u_{0001,0002} = 1$. 

\textit{BRP-PAIRS}'s objective function \eqref{objective3} examines pairs of midshipmen from the same previous company and attempts to minimize the number of student pairs from the same previous company assigned to the same new company. Constraint \eqref{midcombo} forces $u$ to take the value of 1 if both students from the same previous company are reassigned to the same new company, and 0 otherwise.  Constraint \eqref{combobinary} enforces binary restrictions. In addition, constraint \eqref{nocurr} from \textit{BRP-DEV} is also included in \textit{BRP-PAIRS} to ensure that a student does not remain in their previous company.

\section{IMPLEMENTATION}

The three mathematical programming models were implemented in Python 3.8 using the Pyomo 6.0 optimization software package \citep{bynum2021pyomo, hart2011pyomo}. Solutions to the models were obtained with Gurobi 9.1 optimization software \citep{gurobi}, and we used a Windows-10 based computer with two, 2.1 GHz 20 core processors with 192 GB of RAM.  The Gurobi solver used default settings to obtain optimal solutions for the BRP-MIN model.  We used the Gurobi persistent solver interface with Pyomo for the \textit{BRP-DEV} and \textit{BRP-PAIRS} models in order to obtain the details of each feasible solution found.  We stopped the persistent solver when further improvement to the optimality gap plateaued for a significant time.    

For reference purposes, we present Tables \ref{table:Previous-2023} and \ref{table:Previous-2024} with company-level summary statistics of the previous company assignments for the classes of 2023 and 2024.  These tables allow leadership to evaluate how the new reassignment compares to the existing assignment. We note the PRT (Physical Readiness Test) score communicates each student's physical fitness aptitude.

\begin{table}[ht]
  \footnotesize
  \caption{Previous Company-level Summary Statistics: Class of 2023} 
  \centering
  \centering
\begin{tabular}{|c|c|c|c|c|c|} 
\hline
~                    & AOM    & MOM    & \% Male & \% White & PRT score           \\ 
\hline
minimum              & 446.89 & 476.84 & 69.44   & 59.46    & 87.68               \\ 
\hline
maximum              & 679.97 & 648.22 & 78.79   & 80.56    & 92.78               \\ 
\hline
average              & 547.57 & 546.25 & 74.00   & 69.36    & 90.40               \\ 
\hline
standard
  deviation & 61.31  & 39.26  & 2.30    & 5.28     & 1.38                \\ 
\hline
median               & 528.24 & 541.68 & 73.68   & 68.50    & 90.29               \\
\hline
\end{tabular}
  \label{table:Previous-2023}
\end{table}

\begin{table}[ht]
  \footnotesize
  \caption{Previous Company-level Summary Statistics: Class of 2024} 
  \centering
  \centering
\begin{tabular}{|c|c|c|c|c|c|} 
\hline
                     & AOM    & MOM    & \% Male & \% White & PRT score           \\ 
\hline
minimum              & 444.70 & 455.57 & 64.10   & 58.33    & 86.28               \\ 
\hline
maximum              & 664.70 & 633.23 & 74.29   & 82.05    & 92.10               \\ 
\hline
average              & 565.08 & 551.88 & 69.70   & 71.47    & 89.59               \\ 
\hline
standard
  deviation & 56.78  & 43.51  & 1.97    & 5.59     & 1.20                \\ 
\hline
median               & 567.78 & 545.06 & 70.00   & 72.15    & 89.63               \\
\hline
\end{tabular}
  \label{table:Previous-2024}
\end{table}

We designed and implemented the three models near the end of the COVID-19 pandemic in the late spring of 2021 to prepare for the next academic year.  The class of 2021 was not subject to reassignment, as these students will not be in the brigade next year. The class of 2022 was not subjected to reassignment, as USNA leadership determined these students would be in midshipmen leadership roles in the next academic year. Class of 2022 had been sufficiently indoctrinated to the brigade leadership culture prior to the COVID outbreak. Students in the classes of 2023 and 2024 experienced diminished leadership development, as the COVID outbreak affected the brigade from March 2020 until May 2021. Our three models are implemented and separately run for the students of the classes of 2023 and 2024.  USNA leadership chose one of the results and implemented the reassignment for these two classes of students.

\section{RESULTS} \label{results}
Table \ref{table:ModelSize} presents a comparison of each model in terms of size and solution.  First, we note that \textit{BRP-MIN} is smaller than \textit{BRP-DEV} and \textit{BRP-PAIRS} with respect to the number of variables and constraints. \textit{BRP-MIN} obtains a provably optimal solution in three seconds.  We obtain solutions to \textit{BRP-DEV} within 1.5 hours, but the solver was not able to establish a lower bound. Even though we do not obtain a provably optimal solution to \textit{BRP-DEV}, we note the solution metrics are similar to those for the \textit{BRP-MIN} and \textit{BRP-PAIRS} solutions.  Hence, we claim that \textit{BRP-DEV} provides a near-optimal solution. 

\textit{BRP-PAIRS} has the largest number of variables and constraints. Similarly, the solver could not establish a lower bound for \textit{BRP-PAIRS}. However, we claim the solver finds the optimal solution to \textit{BRP-PAIRS} within seven hours. To prove optimality, we first derive \textit{a priori} lower bounds for \textit{BRP-PAIRS} based on observations about the structure of each company (see Appendix). Next, we know that a feasible  solution provides an upper bound on the optimal solution. Therefore, we can claim optimality when the solver obtains a feasible solution with objective function value equal to the lower bound we derive.

\begin{table}[ht]
  \footnotesize
  \caption{Model Size and Solution Comparison} 
  \centering
  \centering
\begin{tabular}{|c|c|c|c|c|} 
\hline
Model~ & Variables & Constraints & Solution Time & Solution Quality  \\ 
\hline
\textit{BRP-MIN}      & 32,911     & 1,984        & 3 seconds     & optimal           \\ 
\hline
\textit{BRP-DEV}      & 33,780     & 4,700        & 1.5 hours     & near-optimal      \\ 
\hline
\textit{BRP-PAIRS}      & 52,474     & 590,000      & 7 hours       & optimal           \\
\hline
\end{tabular}

  \label{table:ModelSize}
\end{table} 
We present the reassignment results from all three models in terms of summary statistics for the class of 2023 and 2024.  We calculate the minimum, maximum, average, median, and standard deviation for a variety of attributes within each newly reassigned company and show the averages for all 30 companies of the brigade.   

Given that \textit{BRP-MIN} minimizes the number of students reassigned to the same previous company, it achieves an optimal objective function value of 0 for both the classes of 2023 and 2024 in 3 seconds. The objective function value of 0 means that no students remained in the same company after reassignment.  

Since \textit{BRP-MIN} achieves the primary goal of eliminating students assigned to the same previous and new companies, the secondary goal of improving the homogeneity of companies comes into focus.  Tables \ref{table:Objective1-2023} and \ref{table:Objective1-2024} show that \textit{BRP-MIN} somewhat improves the homogeneity of each company for the classes of 2023 and 2024, respectively.  However, the final column shows that about 30\%, on average, of students from the same previous company are reassigned to the same new company for both the class of 2023 and 2024. These results are unsurprising because \textit{BRP-MIN} does not control for the number of students moving from a previous company to the same new company.     
 
\begin{table}[ht]
  \footnotesize
  \caption{\textit{BRP-MIN} results for the class of 2023} 
  \centering
  \centering
\begin{tabular}{|c|c|c|c|c|c|c|} 
\hline
~                    & AOM    & MOM    & \% Male & \% White & PRT score & \% from same previous company  \\ 
\hline
minimum              & 489.00 & 495.62 & 66.67   & 60.53    & 87.75     & 19.44                            \\ 
\hline
maximum              & 589.22 & 594.92 & 83.33   & 79.49    & 92.98     & 48.65                            \\ 
\hline
average              & 547.7  & 546.83 & 74.05   & 69.27    & 90.39     & 30.26                            \\ 
\hline
standard deviation   & 30.53  & 29.73  & 4.59    & 5.53     & 1.33      & 7.42                             \\ 
\hline
median               & 554.57 & 547.81 & 75.00   & 69.44    & 90.28     & 28.47                            \\
\hline
\end{tabular}
  \label{table:Objective1-2023}
\end{table}

\begin{table}[ht]
  \footnotesize
  \caption{\textit{BRP-MIN} results for the class of 2024}
  \centering
  \centering
\begin{tabular}{|c|c|c|c|c|c|c|} 
\hline
~                    & AOM    & MOM    & \% Male & \% White & PRT score & \% from same previous company  \\ 
\hline
minimum              & 499.38 & 498.62 & 65.79   & 61.54    & 87.47     & 19.23                            \\ 
\hline
maximum              & 589.97 & 589.1  & 79.49   & 79.49    & 92.23     & 43.42                            \\ 
\hline
average              & 564.6  & 551.61 & 69.72   & 71.51    & 89.59     & 30.85                            \\ 
\hline
standard deviation   & 23.77  & 29.16  & 4.46    & 5.5      & 1.07      & 5.62                             \\ 
\hline
median               & 567.72 & 559.20  & 66.67   & 71.79    & 89.56     & 31.41                            \\
\hline
\end{tabular}
  \label{table:Objective1-2024}
\end{table}

\textit{BRP-DEV} attempts to improve upon the results of \textit{BRP-MIN} by reducing the deviation in the company-level average AOM and MOM scores. Although the solver does not obtain an optimal solution, \textit{BRP-DEV} improves most metrics when compared with \textit{BRP-MIN}.  In Tables \ref{table:Objective2-2023} and \ref{table:Objective2-2024} we observe the standard deviation of the average AOM and MOM scores improved significantly for both the classes of 2023 and 2024.  For the class of 2023, the standard deviation of company AOM scores is reduced from 30 to 12, and MOM score is reduced from 29 to 12.  In the class of 2024, we see the standard deviation of company AOM scores reduced from 23 to 5, and the standard deviation of MOM score reduced from 29 to 5.  The final column shows that about 30\% to 32\% of students from the same previous company are reassigned to the same new company for both the class of 2023 and 2024. Similar to \textit{BRP-MIN}, \textit{BRP-DEV} does not specifically control for the number of students moving from a previous company to the same new company.  

\begin{table}[ht]
  \footnotesize
  \caption{\textit{BRP-DEV} results for the class of 2023} 
  \centering
  \centering
\begin{tabular}{|c|c|c|c|c|c|c|} 
\hline
~                    & AOM    & MOM    & \% Male & \% White & PRT score & \% from same previous company  \\ 
\hline
minimum              & 501.20 &  498.62 &  65.00  & 61.11    & 87.85     & 18.06                            \\ 
\hline
maximum              & 557.53 & 558.44 & 83.78   & 75.68    & 92.85     & 47.22                            \\ 
\hline
average              & 547.93 & 547.03 & 74.06   & 69.29    & 90.38     & 30.91                            \\ 
\hline
standard deviation   & 12.28  & 12.91  & 4.57    & 4.18     & 1.12      & 5.65                             \\ 
\hline
median               & 555.30 & 553.74 & 73.33   & 69.44    & 90.43     & 30.56                            \\
\hline
\end{tabular}
  \label{table:Objective2-2023}
\end{table}

\begin{table}[ht]
  \footnotesize
  \caption{\textit{BRP-DEV} results for the class of 2024}
  \centering
  \centering
\begin{tabular}{|c|c|c|c|c|c|c|} 
\hline
~                    & AOM    & MOM    & \% Male & \% White & PRT score & \% from same previous company  \\ 
\hline
minimum              & 560.05 & 547.18 & 65.79   & 61.54    & 87.89     & 21.79                            \\ 
\hline
maximum              & 578.29 & 564.37 & 79.49   & 79.49    & 91.10     & 43.59                            \\ 
\hline
average              & 564.56 & 551.69 & 69.70   & 71.50    & 89.59     & 32.27                            \\ 
\hline
standard deviation   & 5.46   & 5.21   & 3.25    & 5.05     & 0.82      & 5.02                             \\ 
\hline
median               & 562.12 & 549.47 & 69.23   & 72.74    & 89.54     & 32.05                            \\
\hline
\end{tabular}
  \label{table:Objective2-2024}
\end{table}

\textit{BRP-PAIRS} attempts to address the shortcoming of the first two models, which is the reassignment of multiple students from the same previous company to a new company together.  These students are reassigned to a different numbered company, but they have at least one student from the same previous company with them.  We obtain optimal solutions for \textit{BRP-PAIRS} with objective function values of 227 pairs of students in the same new company for the class of 2023, and 295 pairs for the class of 2024. For \textit{BRP-PAIRS}, we obtain the optimal solution for both the classes of 2023 and 2024 based on the lower bounds derived in the Appendix.

Table \ref{table:Objective3-2023} shows that for the class of 2023, only an average of 10\% of students will find a familiar student from the same previous company reassigned to the new company with them.  Similarly, Table \ref{table:Objective3-2024}  shows that an average of 12\% of students in the class of 2024 will encounter another student from their previous company in their newly assigned company.  \textit{BRP-PAIRS} appears to cut this particular metric significantly when compared to the observed output of the \textit{BRP-MIN} and \textit{BRP-DEV} models.  Unfortunately, \textit{BRP-PAIRS} does not maintain the smaller standard deviation of average academic and military metrics observed in \textit{BRP-DEV}.  

\begin{table}[ht]
  \footnotesize
  \caption{\textit{BRP-PAIRS} results for the class of 2023}
  \centering
  \centering
\begin{tabular}{|c|c|c|c|c|c|c|} 
\hline
~                    & AOM    & MOM    & \% Male & \% White & PRT score & \%  from same previous company  \\ 
\hline
minimum              & 506.30 & 487.47 & 66.67   & 61.11    & 87.46     & 9.72                             \\ 
\hline
maximum              & 584.92 & 593.84 & 83.33   & 77.78    & 93.41     & 12.82                            \\ 
\hline
average              & 547.93 & 546.67 & 73.99   & 69.27    & 90.39      & 10.33                            \\ 
\hline
standard deviation   & 21.92  & 31.76   & 4.29     & 4.52     & 1.32      & 0.80                             \\ 
\hline
median               & 550.32 & 549.80 & 72.97   & 68.93    & 90.63     & 9.72                             \\
\hline
\end{tabular}
  \label{table:Objective3-2023}
\end{table}
\begin{table}[ht]
  \footnotesize
  \caption{\textit{BRP-PAIRS} results for the class of 2024}
  \centering
  \centering
\begin{tabular}{|c|c|c|c|c|c|c|} 
\hline
~                    & AOM    & MOM    & \% Male & \% White & PRT score & \%  from same previous company  \\ 
\hline
minimum              & 504.49 & 498.67 & 66.67   & 60.53    & 86.83     & 11.84                             \\ 
\hline
maximum              & 589.64 & 588.77 & 78.95   & 79.49    & 92.10     & 12.82                            \\ 
\hline
average              & 564.57 & 551.63 & 69.71   & 71.49    & 89.59      & 12.66                            \\ 
\hline
standard deviation   & 21.48  & 25.38   & 3.61   & 6.00     & 1.11      & 0.36                             \\ 
\hline
median               & 565.49 & 552.03 & 68.83   & 71.79    & 89.34     & 12.82                             \\
\hline
\end{tabular}
  \label{table:Objective3-2024}
\end{table}

USNA leadership chose to implement results from \textit{BRP-MIN} to reassign the classes of 2023 and 2024 to new companies, since senior leadership determined that \textit{BRP-MIN} sufficiently achieved all required goals.  The speed (3 seconds) at which the \textit{BRP-MIN} obtains an optimal solution is particularly impressive when compared to the numerous number of man-hours required for a manual reassignment with no guarantee of solution quality. Leadership desired the capability of rapidly solving the BRP to optimality as administrative changes occur to the student membership in some of the subsets in the model.

\section{FUTURE WORK}
The BRP was designed to reassign students to companies as a specific response to the lack of leadership development during the COVID-19 pandemic. However, the BRP model could be modified or expanded to be used for other situations.  We have three suggestions for future use of the BRP.

The BRP could be modified to create an efficient way to perform the initial assignment of students to a company upon first arrival to a service academy. Initial company assignment is made based upon demographics and performance metrics from high school or prior college experience.  No information would exist on military order of merit, except for those students with prior military experience.  Each year, senior leadership completes the initial company assignment for the inbound class by hand.  Applying a mathematical model, such as the ones in this work, would make this process more efficient and successful. The BRP model would have to be changed from a reassignment model to an initial assignment model. All references to the previous company would need to be removed, since new students do not have a previous company assignment.

Another use of the BRP model would be to perform company reassignment on an annual basis. An annual reassignment process could occur at the end of the plebe year for each incoming class of students.  Initial company assignment is made without knowledge of how a student will perform at the collegiate level and information about military performance is usually non-existent.   An annual reassignment of the plebe class at the end of their first academic year when AOM and MOM information is available could be useful. In the academic part of the academy, a similar reassignment of academic advisers occurs.  Plebe students do not declare a major field of study until the second semester of study. A new academic adviser is assigned to each plebe student when an academic major is selected. An annual reassignment could also avoid a high concentration of midshipmen having negative designation (e.g., disciplinary concerns or academic probation). We note that our current models already possess the flexibility to handle such groups in the same way that it disperses a team of varsity sports athletes across all of the companies in the brigade.  

Additional diversity metrics could be added to the model with the goal of further improving the homogeneity of companies. One example area could be geographic diversity. Students attend service academies from all over the United States, and hometown information is available for every student.  Distances could be calculated based upon the zip code of the home address or congressional district of each student.  Distances could be defined either from the school to the student home, or as distances between pairs of students. An objective function could be created to maximize the average distances for students or for pairs of students assigned to a company. The goal of this model would be to create companies that have a mix of students from diverse hometowns.

\section{CONCLUSIONS}
We successfully implemented an IP model for the reassignment of the brigade of midshipmen at USNA. No prior documentation exists for how to perform a reassignment of students to new companies in the brigade, although many senior leaders remember such an event occurring from time to time over the years. In the past, many staff officers reassigned midshipmen to new companies by hand over the course of multiple days.  Alternatively, the mathematical programming models presented in this work can determine an optimal or near-optimal assignment of students in just a few seconds to a few hours depending on the desired outcome. Additionally, all desired requirements for the composition of companies can be enforced through the inclusion of specific constraints.  For example, we enforce upper and/or lower limits on many different metrics, such as company size, academic and military merit score limits, gender and race limits, varsity athlete limits, personal romances, and more.  Brigade leadership no longer has to spend copious amounts of time implementing a sub-par approach to reassigning students to new companies. Our models are shown to be the best possible assignments based on all specified requirements.  

The three different BRP models allow USNA leadership a choice for how to proceed with company-level reassignment.  Each model approaches the question of how to perform the reassignment from different points of view. \textit{BRP-MIN} gives a fast optimal reassignment of students that eliminates the number of students assigned from a previous company to the same numbered company. \textit{BRP-DEV} specifically minimizes the deviation between companies for academic and military order of merit.  \textit{BRP-PAIRS} minimizes the number of students that encounter students in their new company that were also in their old company.  Senior leaders chose fast solution speed provided by \textit{BRP-MIN} for the final implementation, but a compelling case can be made for each of the different models.

\section{ACKNOWLEDGMENTS}
The authors would like to thank the leadership of the USNA Brigade of Midshipmen for using Operations Research analysis to obtain the best possible reassignment of midshipmen to companies.

\section{AUTHOR STATEMENT}
The views expressed in this paper are those of the authors and do not reflect the official policy or position of the United States Navy, Department of Defense (DoD), or the United States Government. The appearance of U.S. DoD visual information does not imply or constitute DoD endorsement.


\bibliographystyle{informs2014} 
\bibliography{main} 

\section{APPENDIX} \label{appendix}
 We derive an \textit{a priori} lower bound for \textit{BRP-PAIRS} which minimizes the number of pairs of students from an old company that are assigned to a new company. We know the previous enrollment for each company for the classes of 2023 and 2024.  Each student in each company must be assigned to a new company different from their previous company, so there are 29 possible new assignments for a student in a previous company.  When the previous company enrollment is greater than 29, it is not possible for each student to be assigned to a new company with no others from the same previous company.  Any enrollment of 30 or more will require some pairs to be assigned to the same new company.  For example, if a previous company has 32 students, at least 32-29 = 3 pairs will need to be assigned to the same new company. 

\begin{center}
\footnotesize
\begin{longtable}{|c|c|c|c|c|c|c|} 
\caption{\textit{BRP-PAIRS} Lower Bounds}\\
\hline
\multicolumn{3}{|c|}{Previous
  Enrollment: Class of 2023} &  & \multicolumn{3}{c|}{Previous Enrollment:  Class of 2024}  \\ 
\hline
Company & Students & Reassignment pairs                    &  & Company & Students & Reassignment pairs                    \\ 
\hline
\endhead
1       & 37       & 8                                     &  & 1       & 39       & 10                                    \\ 
\hline
2       & 35       & 6                                     &  & 2       & 35       & 6                                     \\ 
\hline
3       & 38       & 9                                     &  & 3       & 37       & 8                                     \\ 
\hline
4       & 37       & 8                                     &  & 4       & 38       & 9                                     \\ 
\hline
5       & 38       & 9                                     &  & 5       & 39       & 10                                    \\ 
\hline
6       & 35       & 6                                     &  & 6       & 42       & 13                                    \\ 
\hline
7       & 37       & 8                                     &  & 7       & 36       & 7                                     \\ 
\hline
8       & 38       & 9                                     &  & 8       & 39       & 10                                    \\ 
\hline
9       & 39       & 10                                    &  & 9       & 38       & 9                                     \\ 
\hline
10      & 36       & 7                                     &  & 10      & 40       & 11                                    \\ 
\hline
11      & 39       & 10                                    &  & 11      & 40       & 11                                    \\ 
\hline
12      & 40       & 11                                    &  & 12      & 39       & 10                                    \\ 
\hline
13      & 40       & 11                                    &  & 13      & 40       & 11                                    \\ 
\hline
14      & 36       & 7                                     &  & 14      & 38       & 9                                     \\ 
\hline
15      & 33       & 4                                     &  & 15      & 40       & 11                                    \\ 
\hline
16      & 37       & 8                                     &  & 16      & 39       & 10                                    \\ 
\hline
17      & 38       & 9                                     &  & 17      & 39       & 10                                    \\ 
\hline
18      & 39       & 10                                    &  & 18      & 39       & 10                                    \\ 
\hline
19      & 33       & 4                                     &  & 19      & 39       & 10                                    \\ 
\hline
20      & 36       & 7                                     &  & 20      & 39       & 10                                    \\ 
\hline
21      & 37       & 8                                     &  & 21      & 40       & 11                                    \\ 
\hline
22      & 36       & 7                                     &  & 22      & 39       & 10                                    \\ 
\hline
23      & 36       & 7                                     &  & 23      & 40       & 11                                    \\ 
\hline
24      & 35       & 6                                     &  & 24      & 39       & 10                                    \\ 
\hline
25      & 35       & 6                                     &  & 25      & 38       & 9                                     \\ 
\hline
26      & 37       & 8                                     &  & 26      & 40       & 11                                    \\ 
\hline
27      & 34       & 5                                     &  & 27      & 38       & 9                                     \\ 
\hline
28      & 38       & 9                                     &  & 28      & 39       & 10                                    \\ 
\hline
29      & 33       & 4                                     &  & 29      & 37       & 8                                     \\ 
\hline
30      & 35       & 6                                     &  & 30      & 40       & 11                                    \\ 
\hline
Totals  & 1097     & 227                                   &  & Totals  & 1165     & 295                                   \\
\hline
\end{longtable}
\end{center}

\label{table:Objective3-LB}
Table 10 shows the enrollment of each company before reassignment.  Subtracting 29 from each company enrollment value gives the minimum number of pairs that will need to be assigned to the same new company.  The \textit{BRP-PAIRS} lower bound for the class of 2023 is 227 pairs, and the class of 2024 is 295 pairs.  Unfortunately, even though we can derive a lower bound of \textit{BRP-PAIRS}, we cannot pass the lower bound to the Gurobi solver \citep{najman}. The Gurobi solver is able to obtain the objective function value of 227 for the class of 2023 model in 56 minutes.  Gurobi achieves the objective function value of 295 for the class of 2024 model in 7 hours. Unfortunately, the Gurobi solver does not identify these values as the lower bound.  We use the \textit{BRP-PAIRS} objective function value and the derived lower bound to calculate a relative optimality gap for \textit{BRP-PAIRS}.  We calculate the gap using the standard formula \citep{Sauppe}:
\begin{subequations}
\begin{align}
    & \label{gapformula}
optimality\ gap \ = \ \frac {|best\ solution - best\ bound|}{ best\ bound} \times 100\% \\[10pt]
     & \label{gap2023}
2023 \ optimality\ gap \ = \ \frac {|227 - 227| }{227} \times 100\% = \ 0\%  \\[10pt]
     & \label{gap2024}
2024 \  optimality\ gap \ = \ \frac {|295 - 295 |}{295} \times 100\% = \ 0\%
\end{align}
\end{subequations}

\end{document}